\documentclass[11pt]{article}
\usepackage{amsfonts,amssymb,amsmath,comment}
\usepackage{graphicx,graphics,color}

\usepackage[english]{babel}
\def\R{\mathbb{R}}
\def\Z{\mathbb{Z}}

\def\dd{{\rm d}}

\def\PP{\mathbb{P}}

\def\cD{\mathcal{D}}

\def\cT{\mathcal{T}}
\def\cO{\mathcal{O}}

\newcommand{\Expect}{\mathbf{E}}
\newcommand{\Prob}{\mathbf{P}}

\newlength\mylen
\makeatletter
\@addtoreset{equation}{section}
\makeatother

\newtheorem{proposition}{Proposition}[section]
\newtheorem{theorem}[proposition]{Theorem}

\newcommand{\halmos}{\rule{1ex}{1.4ex}}

\def\R{\mathbb{R}}
\def\Z{\mathbb{Z}}
\def\dd{\mathrm{d}}

\def\cD{\mathcal{D}}

\def\cO{\mathcal{O}}

\begin{document}

\title{Evolution of Discordance}

\author{F.\ den Hollander
\footnote{Mathematical Institute, Leiden University, Einsteinweg 55, 2333 CC Leiden, The Netherlands.}
}

\date{\today}

\maketitle

\begin{abstract}
The present paper is a brief overview of random opinion dynamics on random graphs based on the Ising Lecture given by the author at the World Congress in Probability and Statistics, 12--16 August 2024, Bochum, Germany. The content is a snapshot of an interesting area of research that is developing rapidly.
\end{abstract}



\medskip\noindent
\emph{Acknowledgement:}
The author is supported by the Netherlands Organisation for Scientific Research (NWO) through Gravitation Grant NETWORKS-024.002.003.


\section{Introduction}
\label{s.introduction}

In a friendship network, each individual has an opinion and shares this opinion with friends. Key questions in this context are: How do opinions evolve in a friendship network? In what way does the evolution of the opinions depend on the architecture of the network? What happens when the network itself evolves over time, because friendships are added or removed? What happens when opinions and network are locked in a two-way feedback called \emph{co-evolution}? The answers to these questions are highly challenging.

Opinion dynamics are wel understood on \emph{infinite regular grids}, like the $d$-dimensional Euclidean lattice $\Z^d$ and the $d$-regular tree $\mathcal{T}_d$. On these grids, opinions exhibit a phase transition between \emph{mono-type} equilibria (= consensus) in low dimensions and \emph{multi-type} equilibria (= non-consensus) in high dimensions.

\begin{figure}[htbp]
\begin{center}
\includegraphics[width=0.4\linewidth]{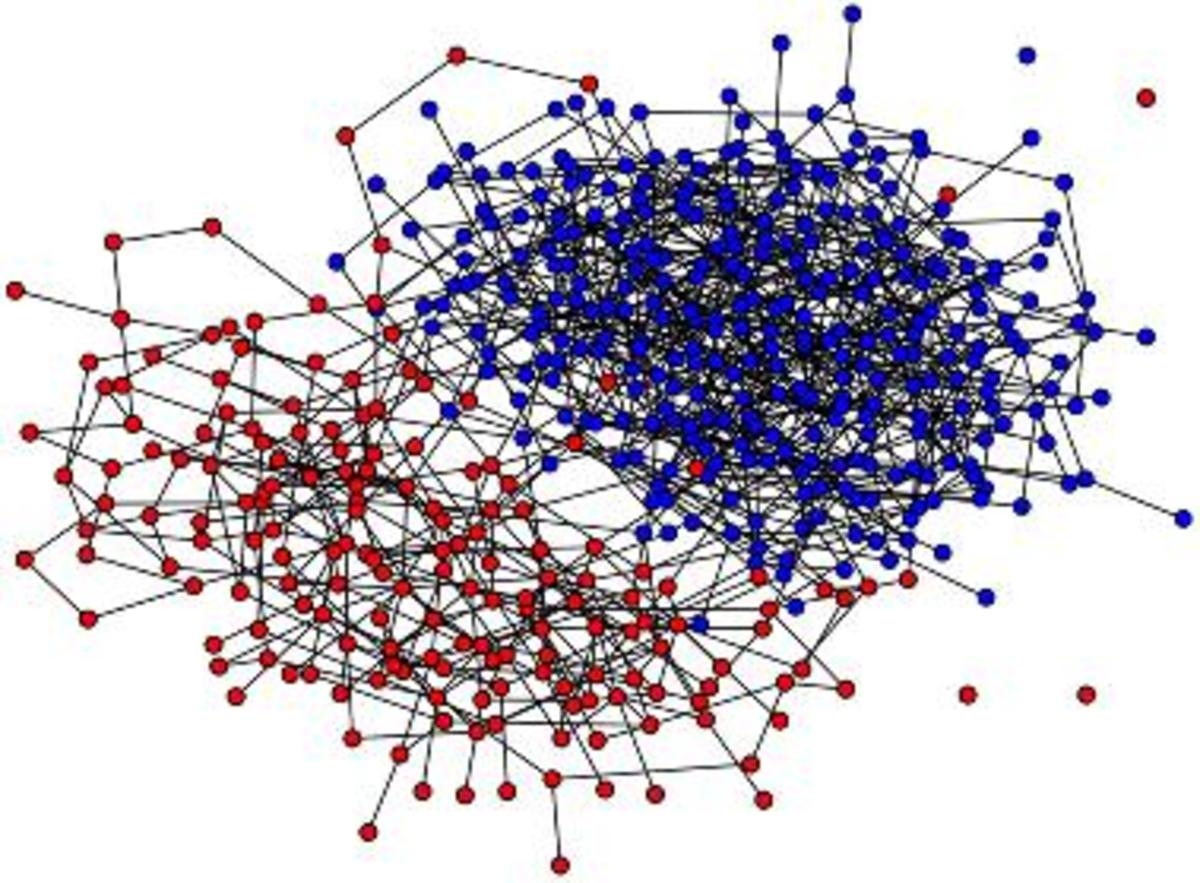}
\end{center}
\vspace{-0.3cm}
\caption{\small Two communities of opinions on a friendship network.}
\label{fig:communities}
\end{figure}

In this lecture we focus on opinion dynamics on a \emph{finite random graph} $G = (V,E)$. We distinguish between \emph{sparse graphs} ($|E| \asymp |V|$) and \emph{dense graphs} ($|E| \asymp |V|^2$). We are interested in the fraction of \emph{discordant edges}, i.e., edges linking vertices with different opinions (see Figure \ref{fig:communities} for an illustration). We will identify \emph{evolution equations} and exhibit \emph{relevant time scales}. 

We further distinguish between two types of dynamics: 
\begin{itemize}
\item[(I)] 
\emph{One-way feedback}. The opinion dynamics depends on the current graph, but the graph dynamics does not depend on the current opinions. In this setting, only consensus is possible, i.e., all individuals eventually have the same opinion.
\item[(II)] 
\emph{Two-way feedback}. Both dynamics depend on each other. In this setting, not only consensus is possible, but also polarisation, i.e., eventually the individuals split into groups having different opinions,
\end{itemize}
Examples of types (I) and (II) are described in Sections~\ref{s.oneway} and \ref{s.twoway}, respectively. 

For background on interacting particle systems we refer the reader to Liggett \cite{L85,L99}, for background on random graphs to van der Hofstad \cite{vdH17,vdH24}. Co-evolution of opinion dynamics and network dynamics poses many mathematical challenges and will be at the \emph{forefront} of network science for years to come.

Ernst Ising was among the first to recognise the importance of \emph{rigorous statistical mechanics} and of attempts to mathematically prove the occurrence of phase transitions in interacting particle systems. A special class is that of \emph{spin systems}, which consist of infinitely many interacting stochastic components that can be in two possible states. Since the pioneering work of Ernst Ising, much progress has been made on understanding spin systems, as highlighted in the contributions elsewhere in this special volume. The classical two-opinion voter model fits perfectly well in this tradition, and is a \emph{gem} of probability theory and mathematical physics. What will be described in Sections~\ref{s.oneway}--\ref{s.twoway} below is a modern-day version of the type of problems that have been driving research on interacting particle systems from the early days onwards.


\section{One-way feedback}
\label{s.oneway}

In this section we focus on examples of models with one-way feedback. Section~\ref{ss.graph} defines the voter model on a graph. Section~\ref{ss.completegraph} looks at the voter model on the complete graph, for which computations can be carried through in full detail. Section~\ref{ss.regulargraph} looks at the voter model on the random regular graph and shows that there are three time scales -- short, modest and long -- exhibiting different behaviour. Section~\ref{ss.directedgraph} shows what happens on directed random graphs, for which certain simplifications can be exploited. Section~\ref{ss.rewiring} investigates what happens on the random regular graph when edges are randomly rewired.


\subsection{Voter model on a graph}
\label{ss.graph}

\begin{figure}[htbp]
\begin{center}
\includegraphics[width=0.4\linewidth]{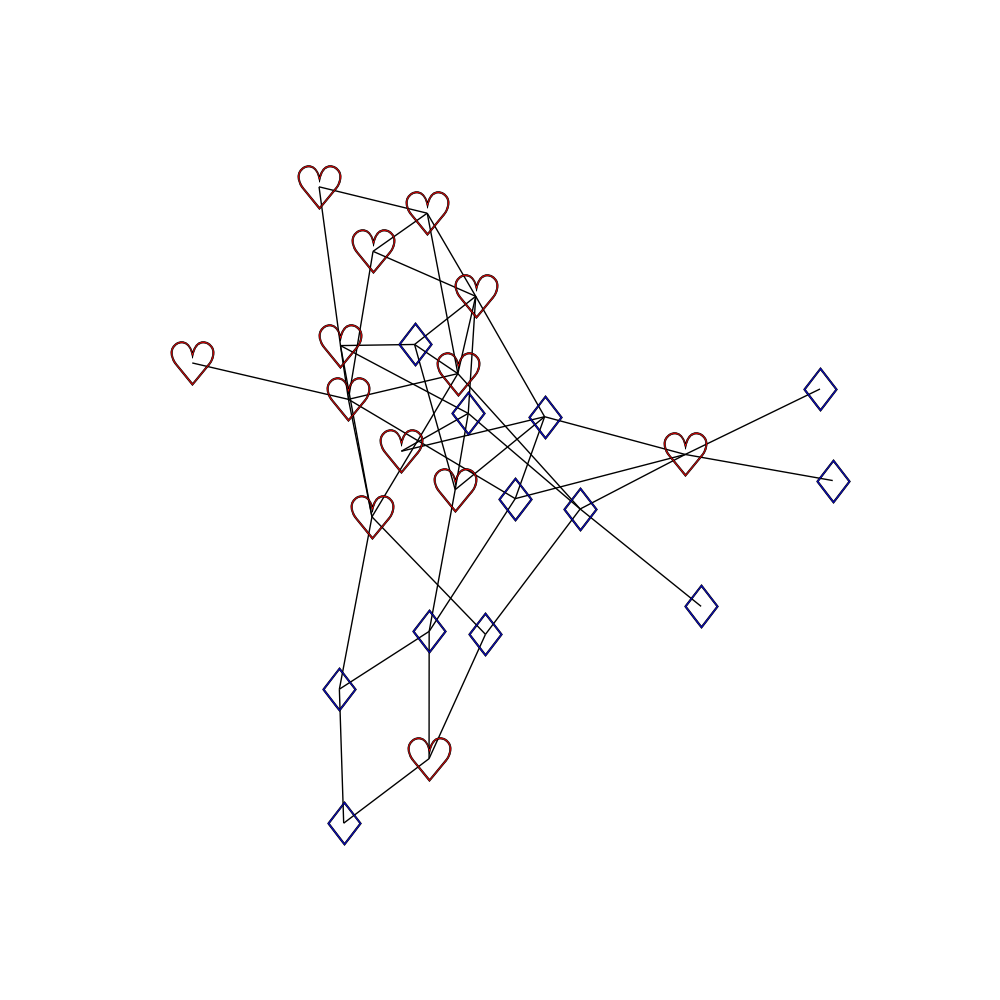}
\end{center}
\vspace{-1.2cm}
\caption{\small A configuration of opinions on a graph. The graph stays fixed, the opinions evolve via random selection and adoption.}
\label{fig:CM}
\end{figure}

The voter model was introduced by Clifford and Sudbury \cite{CS73} and by Holley and Liggett \cite{HL75}. Given a \emph{finite connected} graph $G=(V,E)$, the voter model is the Markov process 
\[
(\xi_t)_{t \geq 0}, \qquad \xi_t = \{\xi_t(i)\colon\,i\in V\},
\]
in which each vertex carries opinion either $\heartsuit$ or $\diamondsuit$, and each vertex at rate $1$ uniformly at random selects a vertex to which it is connected by an edge and adopts its opinion (see Figure~\ref{fig:CM}). The fraction of \emph{discordant edges} is
\[
\cD^N_t = \frac{|D^N_t|}{M}, \qquad D^N_t = \big\{(i,j) \in E\colon\, \xi_t(i) \neq \xi_t(j)\big\},
\]
where $N= |V|$ and $M=|E|$. The latter is an interesting quantity because it monitors the size of the boundary between the two opinions, which is where opinions can change.

The \emph{consensus time} is defined as
\[
\tau = \inf\{t \geq 0\colon\, \xi_t(i)=\xi_t(j)\,\, \forall\, i,j\in V\}.
\]
For finite graphs we know that $\tau<\infty$ with probability $1$, and that consensus is reached at the configurations 
\[
\xi \equiv \heartsuit \quad \mbox{ or } \quad \xi \equiv \diamondsuit.
\] 
The interest lies in determining the \emph{relevant time scales} on which consensus is reached, and \emph{how} it is reached. 

\vspace{-0.75cm}

The voter model is \emph{dual} to a system of coalescing random walks starting from each vertex (see Liggett \cite{L85,L99}). This duality is a major tool in the analysis.


\subsection{Voter model on the complete graph}
\label{ss.completegraph}

\begin{figure}[htbp]
\begin{center}
\includegraphics[width=0.17\linewidth]{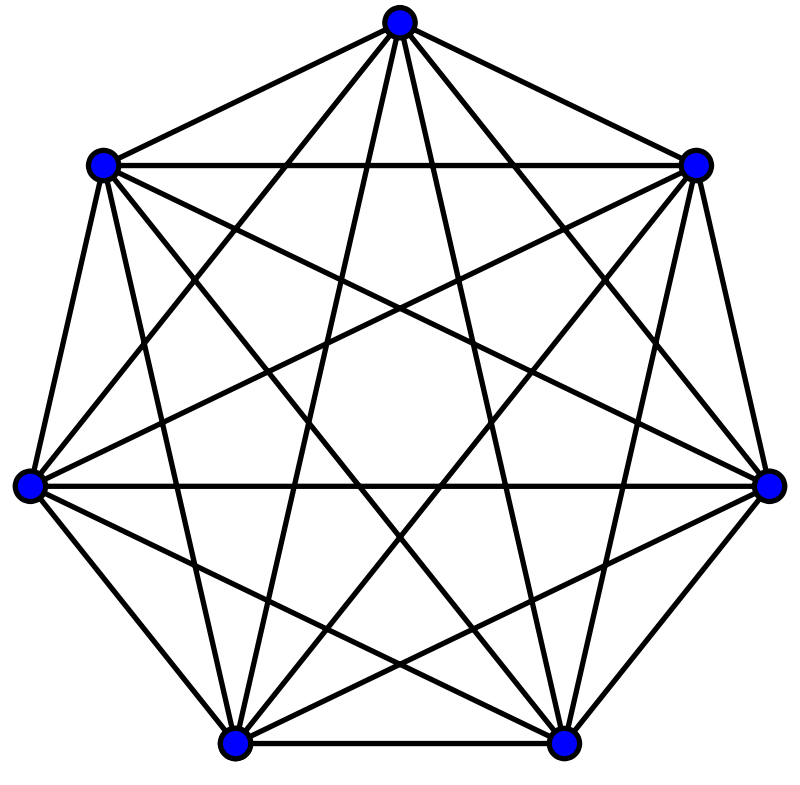}
\end{center}
\vspace{-0.5cm}
\caption{\small The complete graph on $N=7$ vertices.}
\label{fig:completegraph}
\end{figure}

As a prelude we look at the voter model on the \emph{complete graph} (see Figure~\ref{fig:completegraph}), for which computations can be carried out explicitly. Indeed, the number of $\heartsuit$-opinions at time $t$, given by
\[
O^N_t = \sum_{i \in V} 1_{\{\xi_t(i)=\heartsuit\}},
\]
performs a continuous-time random walk on the set $\{0,\ldots,N\}$ with transition rates
\[
\begin{array}{lll}
&n \to n+1 &\text{at rate }  n(N-n) \frac{1}{N-1},\\[0.2cm] 
&n \to n-1 &\text{at rate }  (N-n)n\frac{1}{N-1}.
\end{array}
\]

\begin{figure}[htbp]
\begin{center}
\includegraphics[width=0.35\linewidth]{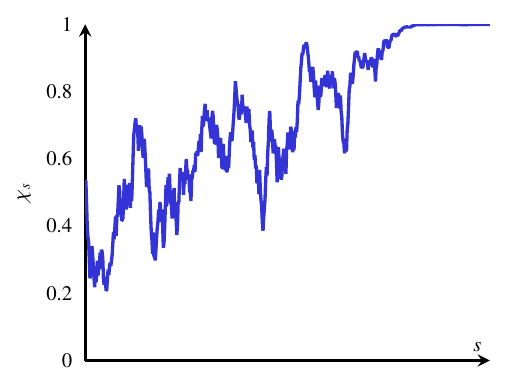}
\end{center}
\vspace{-0.5cm}
\caption{\small A realisation of the Fisher-Wright diffusion, which eventually gets trapped at $0$ or $1$.} 
\label{fig:FW}
\end{figure}

Put $\cO^N_t = \frac{1}{N} O^N_t$ for the fraction of $\heartsuit$-opinions at time $t$. Then it is well-known that
\[
(\cO^N_{sN})_{s \geq 0} \xrightarrow{\rm law} (\chi_s)_{s \geq 0}, \qquad N\to\infty,
\]
where the limiting process is the \emph{Fisher-Wright diffusion}
\[
\dd \chi_s = \sqrt{2 \chi_s (1 - \chi_s)} \,\dd W_s
\]
with $(W_s)_{s \geq 0}$ standard Brownian motion (see Figure~\ref{fig:FW} for a simulation).

\noindent
The number of discordant edges equals
\[
D^N_t = \tfrac12 O^N_t(N-O^N_t).
\]
Put $\cD^N_t=\frac{1}{M}D^N_t$ for the fraction of discordant edges at time $t$, with $M=\binom{N}{2}$ for the complete graph. Since
\[
\cD^N_t = \frac{N}{N-1}\,\cO^N_t(1-\cO^N_t),
\]
it follows that
\[
(\cD^N_{sN})_{s \geq 0} \xrightarrow{\rm law} \big(\chi_s(1-\chi_s)\big)_{s \geq 0}, \qquad N \to \infty.
\]

In the mean-field setting of the complete graph, the fraction of discordant edges is the product of the fractions of the two opinions. The latter property \emph{fails} on non-complete graphs, in particular, on random graphs.


\subsection{Voter model on the random regular graph}
\label{ss.regulargraph}

Consider the \emph{random regular graph} $G_{d,N} = (V,E)$ of degree $d \geq 3$, consisting of 
\[ 
\begin{aligned}
&|V| = N \text{ vertices},\\[0.2cm] 
&|E| = M = \frac{dN}{2} \text{ edges}. 
\end{aligned}
\]

\vspace{-0.5cm}
\begin{figure}[htbp]
\begin{center}
\includegraphics[width=0.2\linewidth]{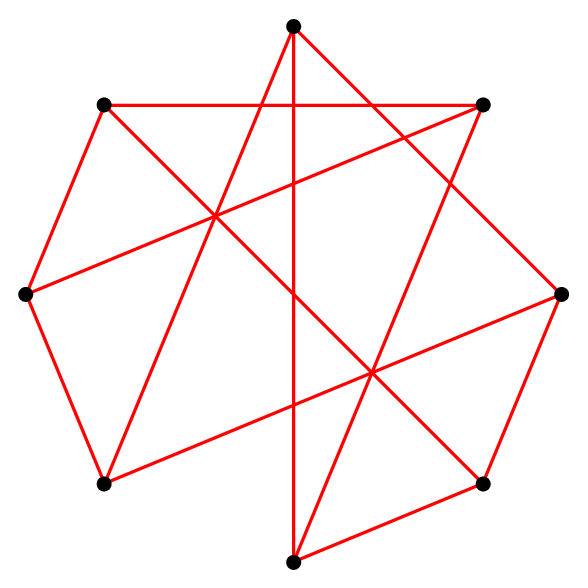} 
\end{center}
\vspace{-0.5cm}
\caption{\small Random regular graph with $d=3$, $N=8$, $M=12$.}
\end{figure}

\noindent
Chen, Choi and. Cox \cite{CCC16} show that
\[
(\cO_{sN}^{N})_{s \geq 0} \xrightarrow{\rm law} (\chi^*_s)_{s \geq 0}, \qquad N \to \infty,
\]
where the limiting process is the Fisher-Wright diffusion
\[
\dd \chi^*_s = \sqrt{2\theta_d\chi^*_s (1-\chi^*_s)}\,\dd W_s
\]
with \emph{diffusion constant}
\[
\theta_d = \frac{d-2}{d-1}.
\]

For $u \in (0,1)$, let ${\bf P}^N_u$ denote the law of $(\cD^N_t)_{t \geq 0}$ when initially the vertices independently have probability $u$ to carry opinion $\heartsuit$. Let $\PP^N$ denote the law of $G_{d,N}$.

\begin{theorem}
\label{thm1}
{\rm (Avena, Baldasso, Hazra, den Hollander, Quattropani \cite{ABHdHQ24}.)}\\
Fix $u\in(0,1)$. Then, for any $t_N \in [0,\infty)$ as $N\to\infty$,
\[
\left| {\bf E}^N_u\left[\cD^N_{t_N}\right] - 2u(1-u)\,f_d(t_N)\,\mathrm{e}^{-2 \theta_d \frac{t_N}{N}}\right| 
\overset{\PP^N}{\longrightarrow} 0
\]
with  
\[
f_d(t) = {\bf P}^{\cT_d}(\tau_{\rm meet} > t),
\]
where ${\bf P}^{\cT_d}$ is the law of two independent random walks on the infinite $d$-regular tree $\mathcal{T}_d$ starting from the endpoints of an edge, and $\tau_{\rm meet}$ denotes their first meeting time. Moreover, for every $\varepsilon>0$, 
\[
\Prob^N_u\left(\left|\cD^N_{t_N}-\Expect^N_u[\cD^N_{t_N}] \right|>\varepsilon \right)
\overset{\PP^N}{\longrightarrow} 0.
\]
\end{theorem}

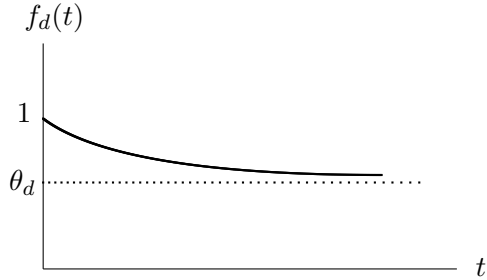
\begin{figure}[htbp]
\begin{center}
\setlength{\unitlength}{0.5cm}
\begin{picture}(10,6)(0,0)
\put(0,0){\line(1,0){11}}
\put(0,0){\line(0,1){6}}
\put(11.5,-0.2){$t$}
\put(-0.5,6.5){$f_d(t)$}
\put(-.7,3.9){$1$}
\put(-.9,2.1){$\theta_d$}
\thicklines\qbezier(0,4)(2,2.5),(9,2.5)
\qbezier[50](0,2.3)(4,2.3)(10,2.3)
\end{picture}
\end{center}
\vspace{-0.5cm}
\caption{\small Plot of $t \mapsto f_d(t)$.}
\label{fig:fdplot}
\end{figure}

\begin{figure}[htbp]
\begin{center}
\includegraphics[width=8cm]{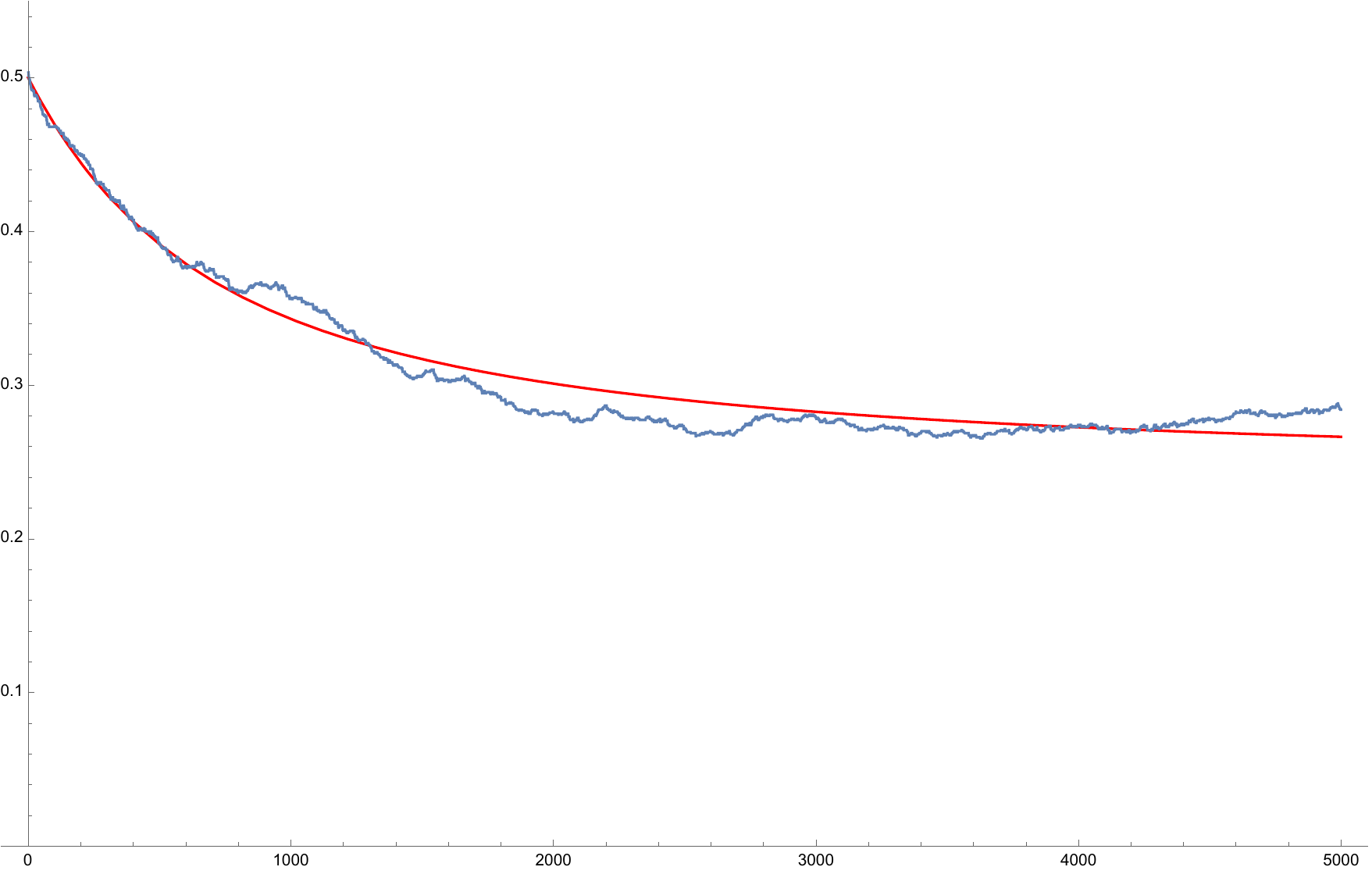}
\end{center}
\caption{\small A single simulation for $N=1000$, $d=3$, $u=0.5$ and $t \in [0,5000]$. In blue the fraction of discordant edges up to $t=5$, in red the function $t \mapsto 2u(1-u)\,f_d(t)$.}
\label{fig:short} 
\end{figure}

\begin{figure}[htbp]
\begin{center}
\includegraphics[width=8cm]{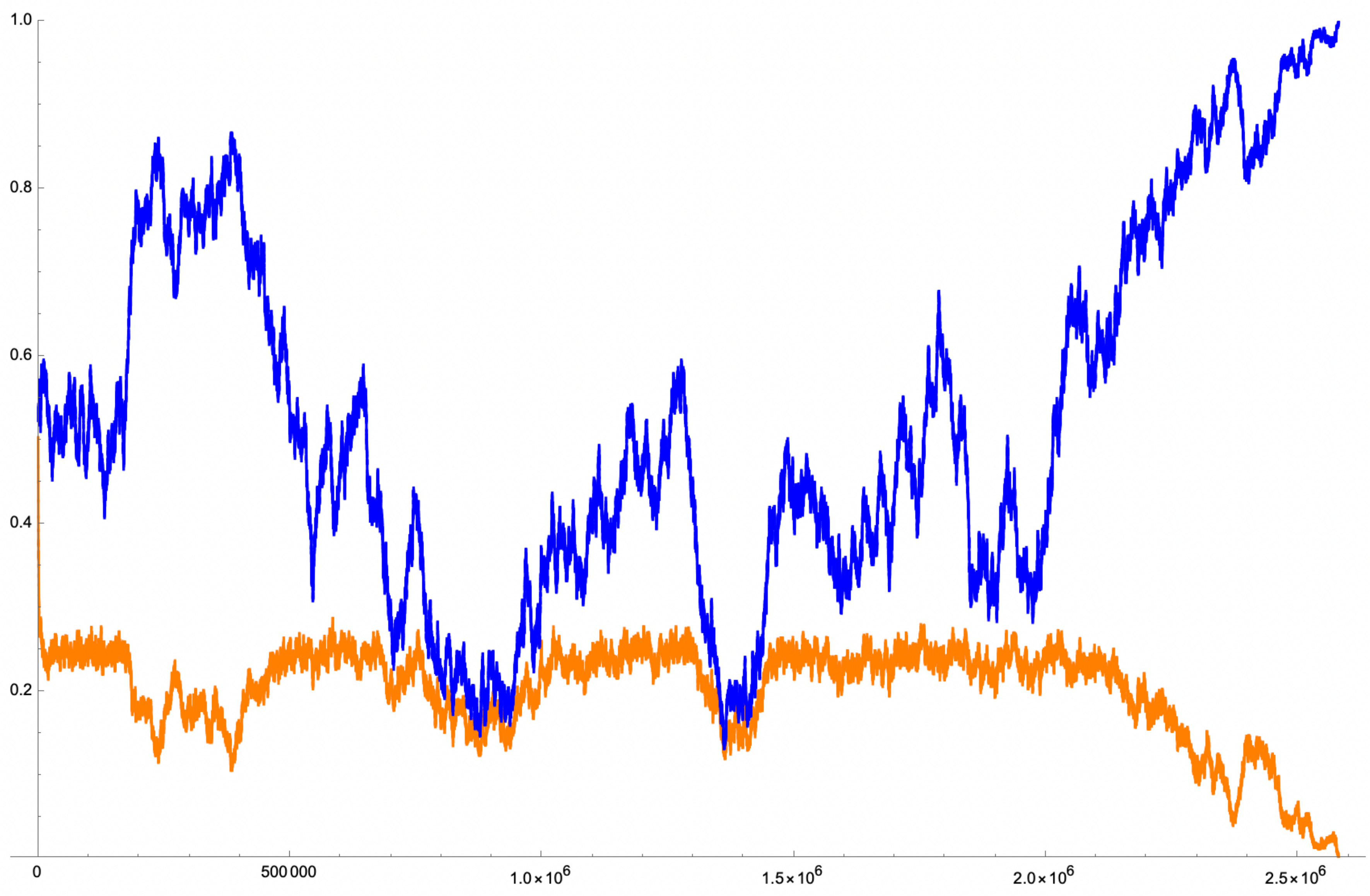}
\end{center}
\caption{\small A single simulation for $N=1000$, $d=3$, $u=0.5$ and $t \in [0,2.5 \times 10^6]$. In blue the fraction of $\heartsuit$-opinions up to consensus, in orange the fraction of discordant edges up to consensus.}
\label{fig:long} 
\end{figure}

\noindent
The function $f_d$ can be computed explicitly and drops from $f_d(0) = 1$ to $f_d(\infty) = \theta_d$ (see Figure \ref{fig:fdplot}). Note that Theorem \ref{thm1} shows that there are \emph{three relevant time scales} (see Figures \ref{fig:short}--\ref{fig:long} for simulations):\\[0.3cm] 
\begin{tabular}{lll}
&short: &$t_N \ll N$,\\
&moderate: &$t_N \asymp N$,\\ 
&long: &$t_N \gg N$.
\end{tabular}

\vspace{0.1cm}\noindent
In a time of order $1$ the fraction of discordant edges drops down from $2u(1-u)$ to $2u(1-u)\theta_d$, where it gets stuck for a time of order $o(N)$. Subsequently, it drops down from $2u(1-u)\theta_d$ to $0$ in a time of order $N$, where it gets stuck after a time of order $\omega(N)$. 

\vspace{1cm}

\begin{theorem}
\label{thm1plus}
{\rm (Avena, Baldasso, Hazra, den Hollander, Quattropani \cite{ABHdHQ24}.)}\\
As $N\to\infty$,
\[
(\cD^N_{sN})_{s \geq 0} \xrightarrow{\rm law} \big(\chi^*_s(1-\chi^*_s)\big)_{s \geq 0}.
\]
\end{theorem}

\noindent
In other words, in the limit as $N\to\infty$ \emph{homogenisation} occurs, i.e., the fraction of discordant edges is again the product of the fractions of the two opinions, just as on the complete graph, despite the fact that the voter model on the random regular graph does \emph{not} fall in the mean-field setting.

Theorems \ref{thm1}--\ref{thm1plus} are examples in a broader program, referred to as the \emph{Cox-Greven finite-systems scheme} (see Cox and Greven \cite{CG90}), where the goal is to understand how interacting particle systems on finite graphs compare with their counterparts on infinite graphs, and to identify what are the relevant time scales for the comparison.  

An early paper listing conjectures for the voter model on random graphs with heavy-tailed degree distributions is Sood, Antal and Redner \cite{SAR08}.


\subsection{Voter model on a directed random graph}
\label{ss.directedgraph}

It remains an open problem to extend Theorems~\ref{thm1}--\ref{thm1plus} to the \emph{configuration model} where the vertex degrees can be different. At present there is not even a conjecture what $\theta_d$ and $f_d(t)$ would become in such a setting. Interestingly, for directed random graphs progress is possible. The setting is the \emph{directed configuration model} with prescribed in-degrees $d^{\mathrm{in}} = (d_i^{\mathrm{in}})_{i=1}^N$ and out-degrees $d^{\mathrm{out}} = (d_i^{\mathrm{out}})_{i=1}^N$, where directed half-edges are paired randomly.

\begin{theorem}
\label{thm2}
{\rm (Avena, Capannoli, Hazra, Quattropani \cite{ACHQ23}, Capannoli \cite{C24pr}.)}\\ 
Under mild conditions on the in-degrees $d^{\mathrm{in}}$ and the out-degrees $d^{\mathrm{out}}$, the same scaling as in Theorem \ref{thm1}--\ref{thm1plus} applies and an explicit formula can be derived both for the profile function $f_{d^{\mathrm{in}},d^{\mathrm{out}}}$ and for the diffusion constant $\theta_{d^{\mathrm{in}},d^{\mathrm{out}}}$.
\end{theorem}

\noindent
For instance, if $d^{\mathrm{in}} = d^{\mathrm{out}}$ (= Eulerian graph), then 
\[
\theta_{d^{\mathrm{in}},d^{\mathrm{out}}} = \left(\frac{m_2}{m_1^2}-1+\sqrt{1-\frac{1}{m_1}}\,\right)^{-1}
\] 
with $m_1,m_2$ the first and the second moment of the limit as $N\to\infty$ of the \emph{empirical distribution of the prescribed degrees}.

It is remarkable that the diffusion constant can be identified for the directed configuration model. Typically directed graphs are harder to deal with than undirected graphs. Apparently, the combination of sparseness and directedness plays out well. It turns out that the dependence of the diffusion constant on $d^{\mathrm{in}}$ and $d^{\mathrm{out}}$ is interesting and reveals that an increase in the volatility of the degrees tends to shorten the consensus time. 


\subsection{Voter model on the random regular graph with random rewiring}
\label{ss.rewiring}

What happens when the random regular graph \emph{itself evolves over time}, for instance, because the edges are randomly rewired? Suppose that every pair of edges swaps endpoints at rate $\nu/2M = \nu/dN$ with $\nu \in (0,\infty)$. Then the rate at which a given edge is involved in a rewiring converges to $\nu$ as $N\to\infty$. The voter model on the resulting dynamic random graph evolves as before: at rate $1$ opinions are adopted along the edges that are currently present.  

\begin{theorem}
\label{thm3}
{\rm (Avena, Baldasso, Hazra, den Hollander, Quattropani \cite{ABHdHQ24pr}.)}\\
The same scaling as in Theorems \ref{thm1}--\ref{thm1plus} applies with the diffusion constant $\theta_d$ in the Fisher-Wright diffusion replaced by
\[
\theta_{d,\nu} = 1-\frac{\Delta_{d,\nu}}{\beta_d},
\]
where $\Delta_{d,\nu}$ is given by the continued-fraction expansion (in the Pringsheim notation)
\[
\Delta_{d,\nu} = \frac{1\,|}{|\,\frac{2+\nu}{\rho_d}} - \frac{1\,|}{|\,\frac{2+2\nu}{\rho_d}} - \frac{1\,|}{|\,\frac{2+3\nu}{\rho_d}} - \dots,
\]
while $\beta_d = \sqrt{d-1}$ and $\rho_d = \frac{2}{d}\sqrt{d-1}$ are constants related to the spectrum of the adjacency matrix of the $d$-regular tree.
\end{theorem}

It is remarkable that the diffusion constant $\theta_{d,\nu}$ can be computed explicitly in terms of a continued-fraction expansion. Apparently, the latter is able to capture the \emph{complexity} inherent in the competition between the dynamics of the opinions and the dynamics of the graph.    

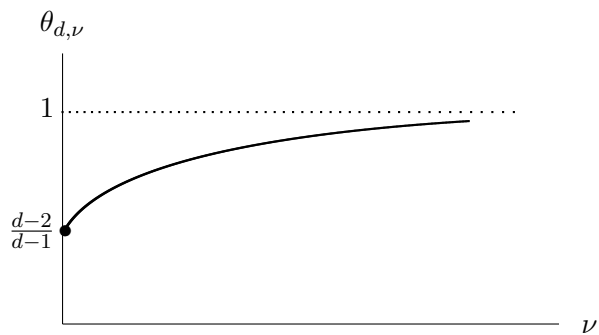
\begin{figure}[htbp]
\begin{center}
\setlength{\unitlength}{0.6cm}
\begin{picture}(10,6)(0,0)
\put(0,0){\line(1,0){11}}
\put(0,0){\line(0,1){6}}
\put(11.5,-0.2){$\nu$}
\put(-0.5,6.5){$\theta_{d,\nu}$}
\put(-.5,4.6){$1$}
\put(-1.2,1.9){$\tfrac{d-2}{d-1}$}
\put(-.1,1.9){$\bullet$}
\thicklines\qbezier(0,2)(1,4),(9,4.5)
\qbezier[50](0,4.7)(4,4.7)(10,4.7)
\end{picture}
\end{center}
\vspace{-0.3cm}
\caption{\small Plot of $\nu \mapsto \theta_{d,\nu}$.}
\label{fig:thetadnuplot}
\end{figure}

Again \emph{homogenisation} occurs, i.e., in the limit as $N\to\infty$ the fraction of discordant edges is the product of the fractions of the two opinions. Moreover, for every $d \geq 3$,
\[
\lim_{\nu \downarrow 0} \theta_{d,\nu} = \frac{d-2}{d-1}, \qquad \lim_{\nu\to \infty} \theta_{d,\nu} = 1,
\]
corresponding to the static random regular graph, respectively, the mean-field random graph. A key observation is that $\nu \mapsto \theta_{d,\nu}$ is \emph{strictly increasing} (see Figure \ref{fig:thetadnuplot}): the rewiring speeds up consensus, even though it does not affect the three relevant time scales: short, moderate, long.   

\vspace{0.3cm}

The identification of the limiting dynamics uses elaborate and delicate coupling arguments between random walks on the random regular graph where edges are randomly rewired and random walks on the regular infinite tree where edges randomly vanish.


\section{Two-way feedback}
\label{s.twoway}

In this section we focus on examples of models with two-way feedback. Section~\ref{ss.coevolution} recalls several attempts to capture co-evolution for opinion dynamics on dynamic random graphs. Section~\ref{ss.double} describes a model for which, in the limit of dense graphs, a \emph{full description} can be given of the co-evolution of the density opinions and the densities of concordant and discordant edges.


\subsection{Attempts at capturing co-evolution}
\label{ss.coevolution}

What happens when the evolutions of the voter model and the graph are locked in a mutual feed-back? Early work for dense graphs, based on simulations and intuitive reasoning, was carried out by Holme and Newman \cite{HN06} and by
Durrett, Gleeson, Lloyd, Mucha, Shi, Sivakoff, Socolar and Varghese \cite{DGLV12}. Depending on the choice of dynamics, a \emph{phase transition} may occur between convergence to a final state where only one opinion survives (=  consensus) or both opinions survive (=  polarisation).

Basu and Sly \cite{BS17} considered the following model proposed in \cite{DGLV12}. Start at time $0$ with the Erd\H{o}s-R\'enyi random graph on $n$ vertices with edge probability $\tfrac12$. For each vertex independently draw from $\{-,+\}$, with probability $\tfrac12$ each, independently of the edges.
\begin{itemize}
\item 
Each edge is assigned a rate-$1$ Poisson clock, independently of everything else. 
\item 
When the clock of an open edge that connects two vertices with different opinions rings, the following happens, with $\beta \in (0,\infty)$ a parameter:
\begin{itemize}
\item 
With probability $n^{-1}\beta$, one of the two vertices is chosen uniformly at random and copies the opinion of the other vertex. 
\item 
With probability $1-n^{-1}\beta$, one of the two vertices is chosen uniformly at random, the edge is disconnected from the vertex that is not selected and reconnected to a vertex chosen uniformly at random from: (I) the set of all vertices (= `rewire-to-random model'); (II) the set of vertices having the same opinion as the selected vertex (= `rewire-to-same model').
\end{itemize}
\end{itemize}
Note that in this model the dynamics ends when there are no edges between vertices with different opinions, which is called the \emph{absorption time} $\tau$. This means that there are two ways in which the dynamics can end: (i) with a positive fraction of vertices holding each opinion (= polarisation); (ii) with all vertices holding the same opinion (= consensus). Classic voter models only allow for case (ii). 

Basu and Sly \cite{BS17} showed that in Models (I) and (II) with high probability, i.e., with probability tending to $1$ as $n \to \infty$, the following dichotomy hold: for $\beta$ small, $\tau$ is of order $n^2$ and polarisation occurs at \emph{equal} fractions of the two opinions, while for $\beta$ large, $\tau$ is of order $n^3$ and polarisation occurs at \emph{different} fractions of the two opinions. As $\beta\to\infty$ the fraction of the minority opinion tends to zero, i.e., there is asymptotic consensus. Moreover, in Model (I) for fixed $\beta$ with high probability the fraction of the minority opinion at time $\tau$ is bounded away from zero. These results were anticipated in \cite{DGLV12} with the help of simulations and approximate computations, and led to the prediction of a \emph{sharp threshold} between the two types of behaviour.

Holme and Newman \cite{HN06} considered a similar dynamics: At each time step choose a vertex uniformly at random.
\begin{itemize}
\item 
With probability $\beta$ choose an edge that is connected to the selected vertex uniformly at random, and reconnect that edge to a vertex that shares the same opinion as the selected vertex chosen uniformly at random.
\item 
With probability $1-\beta$ the selected vertex adopts the opinion of a vertex to which it is connect by an edge chosen uniformly at random. 
\end{itemize}
Holme and Newman \cite{HN06} conjectured a similar phase transition as that conjectured in \cite{DGLV12}, based on simulation. Note that in this model the dynamics also ends when there are no edges between vertices with different opinions. 

Braunsteins, den Hollander and Mandjes \cite{BdHM22} considered a class of graph-valued stochastic processes in which each vertex has a type that fluctuates randomly over time. Collectively, the paths of the vertex types up to a given time determine the probabilities that the edges are active or inactive at that time. The focus is on the evolution of the associated \emph{empirical graphon} in the limit as the number of vertices tends to infinity, in the setting where fluctuations in the graph-valued process are more likely to be caused by fluctuations in the vertex types than by fluctuations in the states of the edges given these types. A \emph{sample-path large deviation principle} is derived, as well as convergence of stochastic processes. The results include a class of stochastic processes where the edge probabilities depend not only on the fluctuations in the vertex types but also on the state of the graph itself. 

Baldassarri, Braunsteins, den Hollander and Mandjes \cite{BBdHM24pr} considered two-opinion voter models on dense dynamic random graphs with the aim to understand the occurrence of consensus versus polarisation. Three models for the joint dynamics of opinions and graphs are studied in which vertices having different opinions are less likely to be connected by an edge than vertices having the same opinion. Functional laws of large numbers for the densities of the two opinions are derived, and \emph{equilibria} are characterised in terms of Beta-distributions.


\subsection{A model where co-evolution can be fully described}
\label{ss.double}

Consider a graph with $N$ vertices.
\begin{itemize}
\item
Each vertex at rate $\eta$ times its current degree rethinks its opinion. When it does so it chooses one of its neighbours uniformly at random and adopts its opinion.
\item
Each pair of vertices at rate $\rho$ rethinks its connection status. When it does so it decides to switch its status as follows:
\[
\begin{tabular}{ll}
&concordant \hspace{0.18cm} \&\,\,$\Bigg\{$\hspace{-0.3cm}
\begin{tabular}{ll} 
connected: &\quad disconnect w.p. $s_{c,1}$,\\[0.2cm] 
disconnected: &\quad connect w.p. $s_{c,0}$,
\end{tabular}
\\
&
\\
&discordant \hspace{0.18cm} \&\,\, $\Bigg\{$\hspace{-0.3cm}
\begin{tabular}{ll}
connected: &\quad disconnect w.p. $s_{d,1}$,\\[0.2cm] 
disconnected: &\quad connect w.p. $s_{d,0}$.
\end{tabular}
\end{tabular}
\]
\end{itemize}
Here, $\eta,\rho \in (0,\infty)$ are parameters that control the speed of the vertex and the edge dynamics, while $s_{c,1},s_{c,0},s_{d,1}s_{d,0} \in (0,1]$ are parameters that control the likelihood for concordant and discordant edges to switch off and on. The limit $N\to\infty$ leads to a limiting dynamics.

\begin{theorem}
\label{thm4}
{\rm (Athreya, den Hollander, R\"ollin \cite{AdHR24}.)}\\ 
Initially, connections and opinions are chosen in an i.i.d.\ manner with densities $p_0, q_0 \in (0,1)$. Let 
\[
\begin{aligned}
p_t &= \text{ density of connected edges at time $t$},\\
q_t &= \text{ density of vertices with opinion $\heartsuit$ at time $t$}.
\end{aligned}
\]
Then
\[
\begin{aligned}
\dd p_t &= \rho\, b(p_t,q_t)\, \dd t,\\[0.2cm]
\dd q_t &= \sqrt{2\eta\, p_t\,q_t(1-q_t)}\, \dd W_t,
\end{aligned}
\]
where $W = (W_t)_{t \geq 0}$ is standard Brownian motion on $\R$, and
\[
b(p,q) = [s_{c,0}\,(1-p) - s_{c,1}\,p]\, [q^2 + (1-q)^2] + [s_{d,0}\,(1-p) - s_{d,1}\,p]\,2q(1-q).
\]
\end{theorem}

\noindent
Thus, the densities of opinions $q_t$ and $1-q_t$ evolve according to a \emph{diffusion}, while the densities of connected and disconnected edges $p_t$ and $1-p_t$ evolve according to a \emph{stochastic flow}, which becomes a \emph{deterministic flow} after consensus is reached, i.e., when $q_t \in \{0,1\}$. The two evolution equations are \emph{coupled}, which is a manifestation of co-evolution. For every $q \in [0,1]$ the drift $p \mapsto b(p,q)$ points to a value $p_c(q) \in (0,1)$, so that eventually consensus occurs. In the boundary case $s_{c,0} = s_{d,0} = 0$ for which disconnected edges cannot connect, the drift points to $0$ and polarisation occurs with a positive probability. 

The identification of the limiting dynamics uses graphon theory and special path topologies. Again, \emph{homogenisation} plays a key role in the background: the opinions rapidly redistribute themselves given the density of the connected and disconnected edges, so that at all times the system is in a \emph{quasi-equilibrium} compatible with the current edge densities. 


\subsection{Simulations}

\begin{figure}[htbp]
\begin{center}
\includegraphics[width=0.5\linewidth]{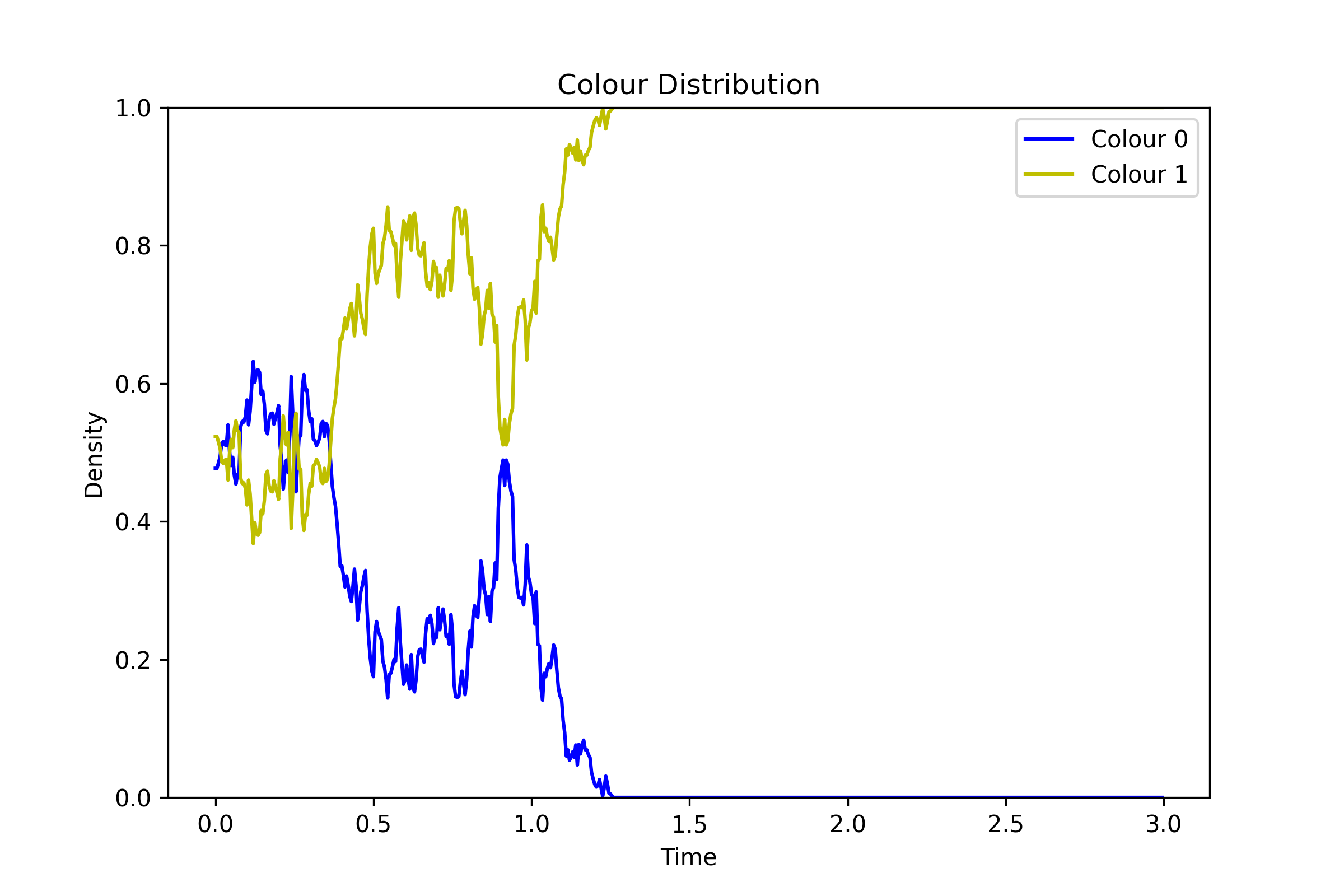}
\end{center}
\vspace{-0.5cm}
\caption{\small Simulations for $N=10^3$ and $\eta=1.0$, $\rho=1.1$ and $s_{c,0} = 1.5$, $s_{c,1} = 0.5$, $s_{d,0} = 0.7$, $s_{d,1} = 2.0$. Plot of the densities of opinions $\diamondsuit$ = colour $0$ (= blue curve); $\heartsuit$ = colour $1$ (= light green curve).}
\label{fig:sim1}
\end{figure}

\begin{figure}[htbp]
\begin{center}
\includegraphics[width=0.475\linewidth]{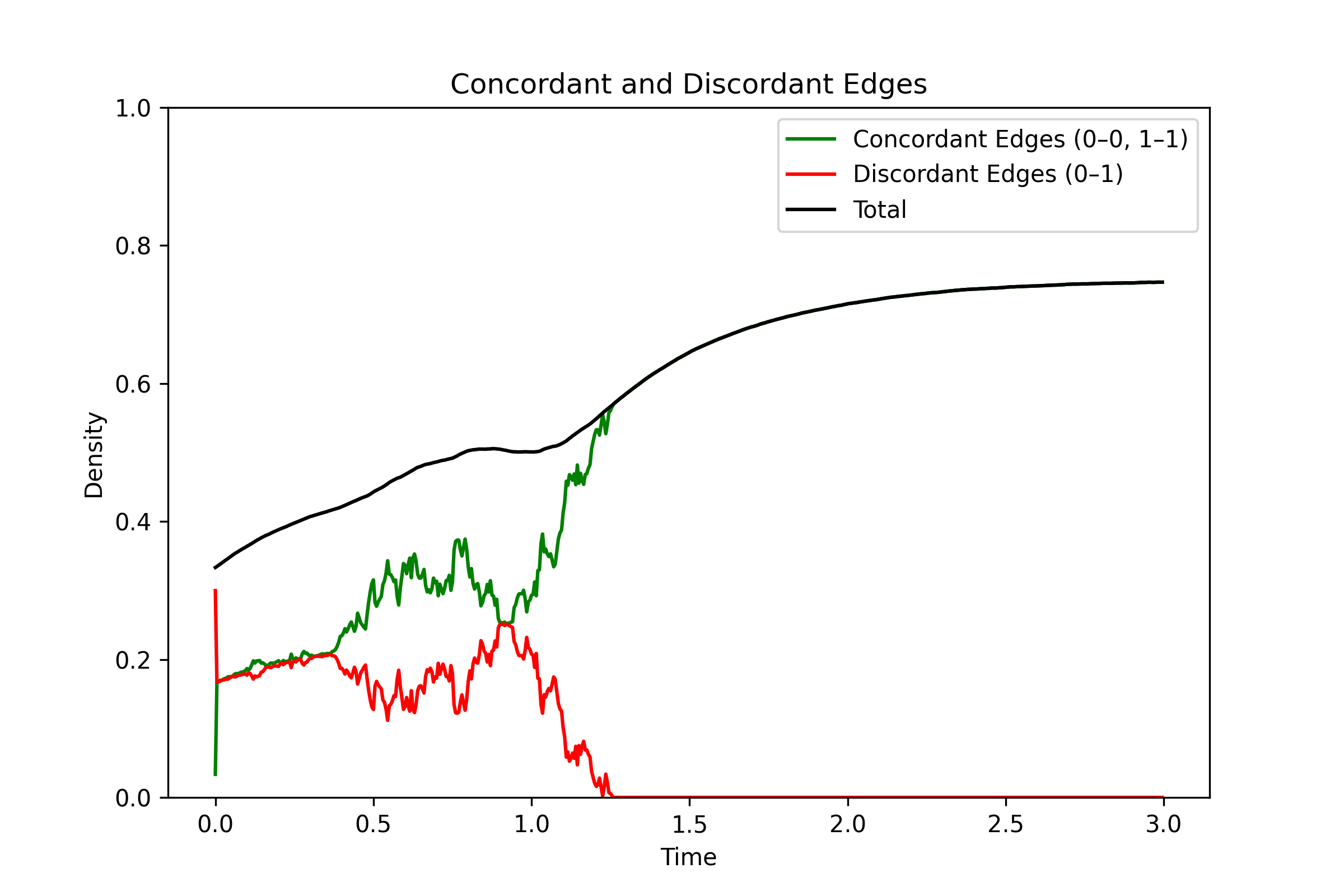}
\includegraphics[width=0.475\linewidth]{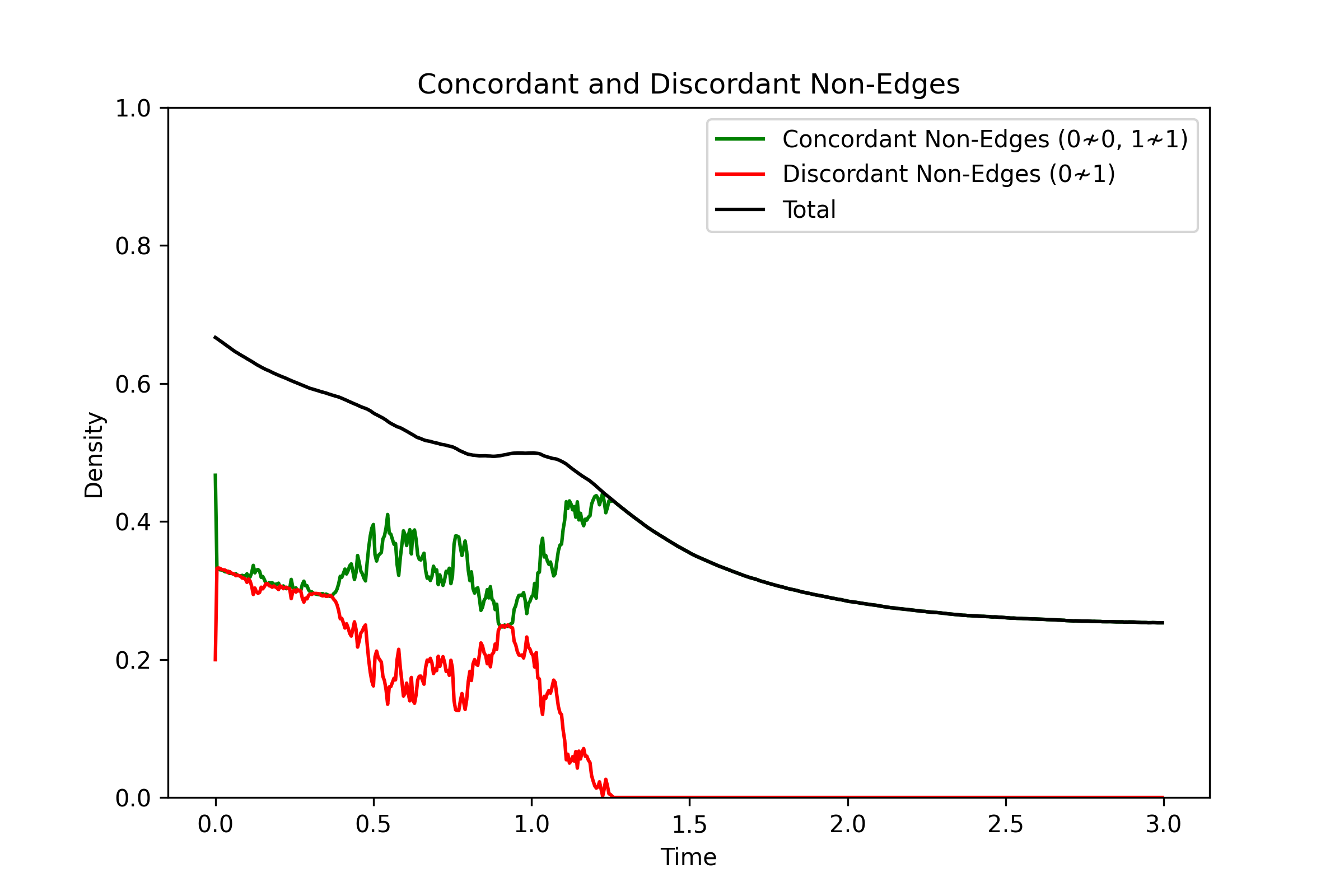}
\end{center}
\vspace{-0.5cm}
\caption{\small \emph{Left}: Plot of the densities of edges (= black curve), concordant edges (= dark green curve), discordant edges (= red curve). \emph{Right}: The same for non-edges.}
\label{fig:sim2}
\end{figure}

Figures~\ref{fig:sim1}--\ref{fig:sim2} show simulations of the pair of evolution equations for $(p_t,q_t)_{t \geq 0}$. The graph is \emph{initialised} as follows. 
\begin{itemize}
\item
Vertices are placed at positions $I_n = \{0,\tfrac{1}{n},\ldots,\tfrac{n-1}{n}\} \subset [0,1]$ and independently receive opinion $0$ or $1$ with probability $\tfrac12$ each. 
\item
Concordant vertices $x,y \in I_n$ are connected with probability $\tfrac{1}{10}\,(1-|x-y|)$, discordant vertices $x,y \in I_n$ are connected with probability $\tfrac{9}{10}\,(1-|x-y|)$.
\end{itemize}
The above choice of parameters is such that discordant edges are more likely to disconnect than to connect, while concordant edges are more likely to connect than to disconnect. The initial graph is therefore \emph{in conflict} with the dynamics: initially discordant edges are more likely to be present than concordant edges. Consequently, there is a \emph{rapid collapse} at the beginning of the evolution. 

\end{document}